\documentclass[onecolumn]{autart}

\usepackage{amsmath,amssymb,amsfonts}
\usepackage{mathtools}
\usepackage{mymath}
\usepackage{MnSymbol}

\usepackage[hidelinks]{hyperref}
\hypersetup{draft=false}

\makeatletter
\let\c@author\relax
\makeatother
\usepackage[backend=bibtex, sorting=none, firstinits=true, abbreviate=true, url=false, isbn=false]{biblatex}
\usepackage{tikz}
\usepackage{pgfplots}
\pgfplotsset{compat=newest}
\usetikzlibrary{plotmarks}
\usepgfplotslibrary{groupplots}

\usetikzlibrary{shapes, arrows, shapes.misc, arrows.meta, positioning, matrix, calc, fit, fadings, patterns}
\usetikzlibrary{calc,patterns,decorations.pathmorphing,decorations.markings,arrows, arrows.meta}
\usepackage{varwidth}
\usepackage[export]{adjustbox}
\usepackage{overpic}

\usepackage{caption}
\usepackage{subcaption}
\usepackage{placeins}

\usepackage{arydshln}
\usepackage{tabto}
\usepackage{siunitx}
\usepackage{nth}

\usepackage{algorithm}
\usepackage[noend]{algpseudocode}
\usepackage{xcolor}

\definecolor{commentgrey}{gray}{0.6}
\newcommand{\algcomment}[1]{%
  \hfill\mbox{}\Comment{\textcolor{commentgrey}{#1}}%
}

\usepackage{todonotes}

\makeatletter
\def\ps@copyright{\let\@mkboth\@gobbletwo
\def\@oddhead{}
\let\@evenhead\@oddhead
\def\@oddfoot{\hfil\thepage\hfil}
\let\@evenfoot\@oddfoot}
\makeatother

\edef\endfrontmatter{%

\unexpanded\expandafter{\endfrontmatter}
\noexpand\endNoHyper
}

\newcommand{\Argmin}{\ensuremath{\mathop{\mathrm{argmin}}}}

\newcommand{\twonorm}[1]{\left\| #1 \right\|_2}
\newcommand{\set}[2]{\left\{ #1 \;\middle|\; #2 \right\}}

\begin{document}
\begin{frontmatter}
\title{Fast-Forwarding Stalling in Dykstra's Algorithm\thanksref{footnoteinfo}}
\thanks[footnoteinfo]{Corresponding author I.~Kempf.}
\author[oxf]{Claudio Vestini}\ead{claudio.vestini@keble.ox.ac.uk}
and
\author[oxf]{Idris Kempf}\ead{idris.kempf@eng.ox.ac.uk}
\address[oxf]{Department of Engineering Science, University of Oxford, Parks Road, Oxford, OX1 3PJ, UK}

\begin{abstract}
Constrained quadratic programs and Euclidean projections are ubiquitous in engineering, arising in machine learning, estimation, control, and signal processing. Dykstra's algorithm is an iterative scheme for computing the Euclidean projection of an initial point onto the intersection of convex sets by successively projecting onto each set. Its low per-iteration computational cost makes it well-suited for solving large-scale or real-time problems where traditional optimisation routines become computationally burdensome. Despite its strong convergence guarantees, Dykstra's algorithm is known to suffer from stalling---arbitrarily long intervals during which the primal iterates remain constant---rendering its runtime unpredictable and severely limiting its applicability in time-critical settings. Focusing on polyhedral constraint sets, we derive a closed-form solution for the length of the stalling period once stalling is detected. This result enables a modified, stall-averse version of Dykstra's algorithm that fast-forwards the stalling period via a single, inexpensive update while preserving convergence guarantees. Numerical experiments demonstrate substantial improvements in convergence behaviour, establishing the proposed method as a practical enhancement for a broad class of projection-based algorithms.
\end{abstract}

\end{frontmatter}

\section{Introduction}

At the core of many problems in optimisation, statistics, and machine learning lies the task of identifying a point that satisfies multiple constraints simultaneously~\cite[Ch. 8]{BOYDCONVEX}. Formally, we consider a Euclidean space $\mathcal{X} \subseteq \Rset^p$ equipped with the Euclidean norm $\twonorm{x}\eqdef (x_1^2+\dots+x_p^2)^{1/2}$, and a finite family of $n$ closed convex subsets $\{\mathcal{H}_i\}_{i=0}^{n-1}\subseteq\mathcal{X}$ with a nonempty intersection $\mathcal{H} := \bigcap_{i=0}^{n-1} \mathcal{H}_i$. Given a point $x^\circ \in \mathcal{X}$, its Euclidean projection $x^\star \reqdef \mathcal{P}_{\mathcal{H}}(x^\circ)$ onto set $\mathcal{H}$ can be found by solving the \textit{constrained quadratic program} (CQP)
\begin{equation} \label{eq:projection}
\begin{array}{ll}
    {\displaystyle \minimise_{x\in\mathcal{X}}}  & \frac{1}{2}\twonorm{x-x^\circ}^2 \\
    \text{subject to} & x \in \bigcap_{i=0}^{n-1}\mathcal{H}_i,
\end{array}
\end{equation}
where $x$ is the decision variable, and $x^\circ$ and sets $\{\mathcal{H}_i\}_{i=0}^{n-1}$ are the problem data. In this paper, we assume $\mathcal{H}$ to be a polyhedral set, meaning $\mathcal{H}=\set{x \in \mathcal{X}}{A x \leq b}$, with $A \in \Rset^{n\times p}$ and $b \in \Rset^n$, and each $\mathcal{H}_i$ a half-space, written as
\begin{align}\label{eq:sets_i}
\mathcal{H}_i\coloneqq\set{x\in\mathcal{X}}{a_i^Tx \leq b_i}, \quad i=0,\dots,n-1,
\end{align}
where vectors $a_i$ are unit normals to each bounding hyperplane $\partial\mathcal{H}_i\coloneqq\set{x\in\mathcal{X}}{a_i^Tx = b_i}$.

The formulation in~\eqref{eq:projection} is widely adopted in optimisation applications, ranging from restricted least squares and isotonic regression in statistics, to signal and image processing (where convex sets encode prior knowledge, like non-negativity or sparsity), control engineering (where system constraints can be modelled using convex sets), and machine learning (where solutions to learning problems must often adhere to structural constraints~\cite{BAUSCHKESWISS}). While in some cases the projection onto $\mathcal{H}$ is easy to compute, for arbitrary polyhedral sets the projection operator $\mathcal{P}_{\mathcal{H}}$ onto the intersection set is intractable and no closed-form solution to~\eqref{eq:projection} is known. In these cases, the solution to~\eqref{eq:projection} can be obtained using standard optimisation software packages. To solve~\eqref{eq:projection}, most solvers require introducing an additional variable $z\eqdef Ax$ (e.g., \cite{osqp} or~\cite{ocpb:16}) that is projected onto the set $\set{z\in\Rset^n}{z \leq b}$, for which an explicit solution exists~\cite{BOYDCONVEX}. However, for large $n$, this variable augmentation can reduce the computational efficiency, particularly in high-dimensional or time-sensitive applications when the projection is part of an overarching algorithm~\cite{OPTIMNESTEROV,KEMPF20206542}.


As an alternative that can be faster in practice~\cite{KEMPF20206542}, \textit{Dykstra's projection algorithm}~\cite{DYKSTRA} is an iterative method designed to solve~\eqref{eq:projection}. It is an extension of the \textit{method of alternating projections} (MAP), which, given an initial iterate, finds a point in the intersection of the half-spaces. Dykstra's algorithm modifies MAP by introducing auxiliary variables to achieve strong convergence to the true projection. Qualitatively, Dykstra's algorithm stores the negative direction vector of projection into auxiliary variables; these are added to the primary variables before subsequent projections and, as a result, the primary iterates of Dykstra's method asymptotically converge to the optimal solution of~\eqref{eq:projection}. Although, as in standard optimisation packages, Dykstra's method introduces additional variables, these appear only in vector additions rather than in matrix-vector multiplications. This reduces the computational load and can make Dykstra's method a highly valuable algorithm to deploy in real-time, high-throughput applications for finding fast (approximate) solutions to CQPs. However, despite its strong theoretical guarantees, Dykstra's algorithm is known to suffer from a \emph{stalling} problem~\cite{DYKSTRASTALLING}. In such cases, the iterates of Dykstra's method remain unchanged for a number of iterations. The duration of the stalling period is not known \emph{a priori}, and, given the choice of starting point, can be made arbitrarily long. In general, the algorithm cannot be run for a fixed number of iterations with a guarantee that the output will be closer to the projection than the initial point, which can hinder deployment of Dykstra's method in practical applications.

In this paper, we address the stalling problem of Dykstra's method for polyhedral sets. Exploiting simplifications that arise in the polyhedral case, we derive the exact duration of the stalling period once stalling is detected (Theorem~\ref{thm:nstall}). This quantity enables a modified version of Dykstra's algorithm (Algorithm~\ref{alg:fastforward}) that fast-forwards through the stalling period, thereby improving convergence properties. Crucially, this modification preserves all convergence guarantees of the original algorithm (Corollary~\ref{cor:convergence}). This development represents a step towards practical deployment of Dykstra's algorithm in real-time settings, where its low-overhead projection could offer significant performance advantages. The structure of the paper is as follows. Section~\ref{sec:background} reviews MAP, Dykstra's algorithm, and narrows the analysis to polyhedral sets, introducing the stalling phenomenon. Section~\ref{sec:main_result} proposes the fast-forward modification, and Section~\ref{sec:example} presents supporting numerical experiments, which demonstrate the effectiveness of the solution. Finally, Section~\ref{sec:conclusion} concludes with a discussion of potential applications, extensions, and refinements.

\emph{Notation}.
We denote the transpose of a vector $a$ as $a^T$, and its $\ell_{2}$-norm as $\twonorm{a}$. The interior of a set $\mathcal{X}$ is denoted as $\text{int}\,\mathcal{X}$, and its boundary is denoted as $\partial\mathcal{X}$. Function composition is written as $g \circ h$, and the modulo and ceiling operators are denoted by $[\cdot]$ and $\lceil \cdot \rceil$, respectively. The absolute iteration index is denoted as $m$, the half-space index (i.e. which half-space we are currently considering) is denoted as $i$, while the half-space count (i.e. how many half-spaces there are in total) is denoted as $n$, and is fixed. We name an update with an increase by $1$ in $m$ ($m \leftarrow m+1$) an \emph{iteration}, and refer to an advancement of $n$ iterations ($m \leftarrow m+n)$ as a \emph{cycle}.

\section{Background} \label{sec:background}

The MAP algorithm is an iterative method to find a point $x_m$ in the intersection region, $x_m  \in\bigcap_{i=0}^{n-1}\mathcal{H}_i$, by cyclically projecting onto the convex sets $\mathcal{H}_i$~\cite{NEUMANN,BREGMAN1965}. Starting from an arbitrary point $x_0=x^\circ$, MAP generates a sequence of iterates $\{x_m\}$ by computing successive projections onto each of the half-spaces,
\begin{align}\label{eq:map}
    x_{m+1} = (P_{\mathcal{H}_{n-1}} \circ \dots \circ P_{\mathcal{H}_0}) (x_m).
\end{align}
When each $\mathcal{H}_i$ is an affine subspace, it can be shown that the MAP also solves problem~\eqref{eq:projection}~\cite{NEUMANN}. For general closed convex sets, however, the sequence generated by successive iterations of~\eqref{eq:map} is only guaranteed to converge to a point within the intersection $\mathcal{H}$, and not to the Euclidean projection $x^\star$~\cite{BREGMAN1965}.

Dykstra's method extends~\eqref{eq:map} by introducing auxiliary variables $e_m \in \mathcal{X}$ that store the negative error vectors associated with the projection of a primary variable $x_m$ onto a given half-space. The algorithm proceeds by generating a series of primary iterates \{$x_{m}$\} and auxiliary iterates \{$e_{m}$\} using the following scheme:
\begin{subequations}\label{eq:dykstra}
\begin{align}
x_{m+1}&=\mathcal{P}_{\mathcal{H}_{[m]}}\left(x_{m}+e_{m-n}\right),\label{eq:dykstra:proj}\\
e_m&=e_{m-n}+x_{m}-x_{m+1}\label{eq:dykstra:error},
\end{align}
\end{subequations}
where $[m]$ represents the modulus operator $[m]\eqdef m \, \text{mod} \, n$. The decision variable is initialised as $x_0=x^\circ$ and the auxiliary variables $e_m$ as $e_{-n} = e_{-(n-1)} = \dots = e_{-1} = 0$.
The Boyle-Dykstra theorem~\cite{DYKSTRA,DYKSTRAPERKINS} implies asymptotic convergence to the projection, so that $\lim_{m\rightarrow\infty}\anynorm{x_m-\mathcal{P}_\mathcal{H}(x^\circ)}=0$. For any finite number of iterations $m$, there is no guarantee that the iterates $x_m$ lie in the intersection $\mathcal{H}$, nor that $x_m\neq x^\circ$. Note that the MAP update~\eqref{eq:map} can be recovered from Dykstra's update~\eqref{eq:dykstra} by setting the auxiliary variables to zero in~\eqref{eq:dykstra:error}, i.e.\ $e_m=0$ for all $m$.

For polyhedral sets as defined in~\eqref{eq:sets_i}, the projection step~\eqref{eq:dykstra:proj} of Dykstra's method can be simplified to
\begin{align}\label{eq:dykstra:proj:poly}
x_{m+1}=
\begin{cases}
x_{m}+e_{m-n} & \text{if } x_{m}+e_{m-n}\in\mathcal{H}_{[m]}\\
x_{m} - \left(x_{m}^T a_{[m]} - b_{[m]}\right) a_{[m]} & \text{otherwise}
\end{cases},
\end{align}
and the update for the auxiliary vector~\eqref{eq:dykstra:error} to
\begin{align}\label{eq:dykstra:error:poly}
e_m=
\begin{cases}
0 & \text{if } x_{m}+e_{m-n}\in\mathcal{H}_{[m]}\\
e_{m-n}+\left(x_{m}^T a_{[m]} - b_{[m]}\right) a_{[m]} & \text{otherwise}
\end{cases}.
\end{align}
The auxiliary vector $e_m$ is either zero or parallel to the half-space norm $a_{[m]}$, so it can be characterised via a scalar quantity $k_m \in \mathbb{R}$ as $e_m = k_m a_{[m]}$ with 
\begin{align}\label{eq:km}
k_m = 
\begin{cases}
0 & \text{if } x_{m}+k_{m-n}a_{[m]}\in\mathcal{H}_{[m]}\\
k_{m-n} + x_m^T a_{[m]} - b_{[m]} & \text{otherwise}
\end{cases}.
\end{align}
The convergence of Dykstra's iterates to the Eucledian projection has been analysed in~\cite{DYKSTRAPERKINS,DYKSTRAPOLY,DYKSTRAPOLY2} for polyhedral sets. The analysis in~\cite{DYKSTRAPOLY} is based on partitioning the collection $\{\mathcal{H}_i\}_{i=0}^{n-1}$ into an \emph{inactive} subset---half-spaces containing the projection $x^\star$ in their interior $\text{int}\,\mathcal{H}_i$---and an \emph{active} subset---half-spaces containing the projection on their boundary $\partial\mathcal{H}_i$. The two subsets, $\mathcal{A}$ and $\mathcal{B}$, respectively, are defined as:
\begin{align}
&\mathcal{A} \eqdef \set{i\in\lbrace 0,\dots,n-1\rbrace}{x_\infty\in \partial\mathcal{H}_i},
&\mathcal{B}\eqdef \lbrace 0,\dots,n-1\rbrace\backslash \mathcal{A}=\set{i\in\lbrace 0,\dots,n-1\rbrace}{x_\infty\in \text{int}\,\mathcal{H}_i},
\end{align}
where $x_\infty=\lim_{m\rightarrow\infty}x_m$. It can be shown that there exists a number $N_1$ such that when $[m]\in \mathcal{B}$ at iteration $m\geq N_1$, it follows that $x_m=x_{m-1}$ and $e_m=0$, i.e.\ the half-spaces that become inactive remain so~\cite[Lemma~3.1]{DYKSTRAPOLY}. Furthermore, there exists $N_2\geq N_1$ such that whenever $n\geq N_2$, it holds that $\twonorm{x_{m+n}-x_\infty}\leq\alpha_{[m]}\twonorm{x_m-x_\infty}$, where $0\leq\alpha_{[m]}<1$ are scalars related to angles between half-spaces~\cite[Lemma~3.7]{DYKSTRAPOLY}. The integer $N_2$ describes the iteration from which on the algorithm has determined the inactive half-spaces. It then follows that there exist constants $0\leq c < 1$ and $\rho > 0$ such that $\twonorm{x_m -x_\infty} \leq \rho\, c^m$~\cite[Thm.~3.8]{DYKSTRAPOLY}. The constant $c$ can be estimated from the smallest $\alpha_{[m]}$, which is characterised by the angle between certain subspaces (subspaces formed by active half-spaces). The constant $\rho$, however, depends on an unknown $N_3\geq N_2$ and on $x^\circ$, and can therefore not be computed in advance~\cite{DYKSTRAPERKINS,XUPOLY}. 

Although the Boyle-Dykstra theorem implies guaranteed asymptotic convergence to the optimal solution $x^\star$, the iterates in Dykstra's algorithm are not guaranteed to converge monotonically. It is in fact possible to establish initial conditions for which the iterates exhibit a cyclic repetition pattern, known as \emph{stalling}.
As a motivating example, in~\cite{DYKSTRASTALLING} the behaviour of Dykstra's method is analysed for two polyhedral sets in $\mathbb{R}^2$, namely an origin-centred box intersected with a line. For this specific example, the authors give conditions on Dykstra's algorithm for (i) finite convergence, (ii) infinite convergence, and (iii) stalling followed by infinite convergence, identifying sets of initial conditions to achieve each of the convergence mechanisms.

\section{Main result}
\label{sec:main_result}

Based on the discussion in~\cite{DYKSTRASTALLING}, we define stalling as follows.

\begin{defn}[Stalling]\label{def:stalling}
Given $m\geq n-1$, a stalling period is defined as those $i\in\lbrace 0,\dots,m_{stall}\rbrace$ with $m_{stall}\geq n$ for which $x_{m+i}=x_{m+i-n}$.
\end{defn}

During stalling, the primary variables $x_{m+i}$ remain constant, and only the auxiliary variables $e_{m+i}$ change: these are modified by constant increments $x_m-x_{m+1}$. The stalling period therefore continues until the auxiliary variable $e_{m+i}$ becomes such that $x_{m+i}+e_{m-n+i}\in\mathcal{H}_{[m+i]}$ for some $[m]$, i.e. the half-space associated with stalling becomes inactive. For polyhedral sets, the length of the stalling period can be pre-computed as soon as the algorithm encounters stalling, using the half-spaces that are active at this iteration. The stalling period continues until one of these half-spaces becomes inactive, at which point stalling terminates. This is formalised in Theorem~\ref{thm:nstall}.

\begin{thm}[Length of Stalling Period]\label{thm:nstall}
Suppose that stalling starts at iteration $m$ as in Definition~\ref{def:stalling} and denote the active half-spaces by $\mathcal{A}\coloneqq\set{i\in\lbrace 0,\dots,n-1\rbrace}{x_{m+i}+e_{m-n+i}\not\in\mathcal{H}_{[m+i]}}$. The length of the stalling period is
\begin{align}\label{eq:nstall}
N_{stall}\coloneqq\min_{i\in\mathcal{S}} \left\lceil \frac{k_{m-n+i}}{b_{[m+i]}-x_{m+i}^T a_{[m+i]}}\right\rceil,
\end{align}
where $\lceil\cdot\rceil$ is the ceiling operator, $k_{m-n+i}=e_{m-n+i}^T a_{[m+i]}$, and $\mathcal{S}\coloneqq\set{i\in \mathcal{A}}{x_{m-n+i}^T a_{[m+i]} - b_{[m+i]} < 0}$.
\end{thm}

\begin{pf}
Note that the existence of a finite number $N_{stall}$ is a consequence of the Boyle-Dykstra theorem. Following the condition in~\eqref{eq:dykstra:proj:poly}, the stalling period will terminate once an active half-space becomes inactive, which we can express as $x_{m+i}+e_{m-n+i}\in\mathcal{H}_{[m+i]}$ for some half-space index $i\geq n$ (Definition~\ref{def:stalling}). As a consequence of the scalar update in~\eqref{eq:km}, the only half-spaces that can become inactive are those for which the quantity $x_{m+i}^T a_{[m+i]} - b_{[m+i]}$ is negative. By Definition~\ref{def:stalling}, this quantity is constant during stalling. The length of the stalling period can therefore be obtained from choosing the smallest possible integer $N_i$ for which
\begin{align*}
k_{m-n+i}+N_i\left(x_{m+i}^T a_{[m+i]} - b_{[m+i]}\right) < 0,\qquad i\in\mathcal{A},
\end{align*}
so that $N_i=\left\lceil k_{m-n+i} / \left(b_{[m+i]}-x_{m+i}^T a_{[m+i]}\right)\right\rceil$, $N_{stall}=\min_{i\in\mathcal{S}}N_i$, and $i_{stall}=\Argmin_{i\in\mathcal{S}} N_i$. If multiple $i_1<\dots<i_j$ are such that $N_{i_1}=\dots = N_{i_j}$, then according to~\eqref{eq:dykstra}, half-space at index $i_{stall}:= i_1$ is discarded first.\qed
\end{pf}
Note that $N_{stall}$ in Theorem~\ref{thm:nstall} refers to the number of $n$-step cycles of~\eqref{eq:dykstra} to skip. Following the end of the stalling period, one half-space is discarded from the active half-space and $\mathcal{A}$ redefined, upon which the algorithm may again encounter a stalling condition. By Theorem~\ref{thm:nstall}, the stalling period can be fast-forwarded by applying constant increments to the auxiliary variables $k_{m}$. Theorem~\ref{thm:nstall} and Dykstra's algorithm are combined in Algorithm~\ref{alg:fastforward}.

\algnewcommand\Inputs{\item[\textbf{Inputs:}]}
\algnewcommand\Outputs{\item[\textbf{Output:}]}

\begin{algorithm}
\caption{Dykstra's projection algorithm for polyhedral sets with fast forwarding.}
\label{alg:fastforward}
\begin{algorithmic}[1]
    \Inputs $x^\circ, \{a_i, b_i\}_{i=0}^{n-1}, N_{\max}, \epsilon_{\text{stall}}$
    \Outputs projection $x_m$
        \State $x_0 \leftarrow x^\circ$; $e_j \leftarrow 0\,\forall j\in\lbrace -n,\dots,-1\rbrace$;  $m \leftarrow 0$ 
        \While{$m < N_{\max}$}
            \State $x_{m+1} \leftarrow \mathcal{P}_{[m]}(x_m + e_{[m-n]})$
            \State $e_{m} \leftarrow e_{m-n} + x_m - x_{m+1}$
            \If{$m \ge n-1 \;\textbf{and}\;\twonorm{x_{m+1-i} - x_{m+1-i-n}} < \epsilon_{\text{stall}}\;\forall i\in\{0,\dots,n-1\}$}\label{alg:start} \algcomment{Stalling detected per Definition~\ref{def:stalling}}
                \State Obtain $N_{stall}$ and $i_{stall}$ from Theorem~\ref{thm:nstall}
                \For{$j \in \lbrace 0,\dots,n-1\rbrace$}
                    \State $N_{skip}\leftarrow N_{stall}$ \textbf{if} $[m-j]\leq i_{stall}$ \textbf{else} $N_{stall}-1$
                    \State $e_{m-j} \leftarrow e_{m-j} + N_{skip}(a_{[m-j]}^T x_{m-j} - b_{[m-j]})a_{[m-j]}$\algcomment{Fast-forward auxiliary variable}
                \EndFor
                \State $m \leftarrow m - i_{stall}$\label{alg:end}\algcomment{Adjust iteration counter}
            \Else
                \State $m \leftarrow m + 1$
            \EndIf
        \EndWhile
        \State \Return $x_{m-1}$
\end{algorithmic}
\end{algorithm}

The algorithm proceeds identically to the standard Dykstra update until it detects stalling by verifying the stationarity condition $x_m = x_{m-n}$ for all half-spaces, as per Definition~\ref{def:stalling}. Once stalling is detected at iteration $m$, the algorithm executes its fast-forwarding phase: after performing the standard variable update for the current iteration, it computes the exact number of required stalling cycles, $N_{stall}$, and the half-space at which stalling ends, $i_{stall}$, by calling a function that implements Theorem~\ref{thm:nstall}, determining the number of full $n$-step cycles to skip. For each of the $n$ auxiliary variables $e_j$, the algorithm then computes the constant increment $\Delta e_j = (a_j^T x_{\text{stall}} - b_j)a_j$ and applies a single update: $e_j \leftarrow e_j + N_{\text{skip}} \Delta e_j$, where $N_{\text{skip}}=N_{\text{stall}}$ if $i_{\text{stall}}\geq j$ and $N_{\text{skip}}=N_{\text{stall}}-1$ otherwise. This ``jumps'' the dual variables forward to their state following the stalling period. Finally, the main iteration counter $m$ is adjusted. For practical purposes, in this implementation Algorithm~\ref{alg:fastforward} is executed for a fixed number of iterations $N_{\text{max}}$, the stationarity condition is enforced with a small nonzero tolerance $\epsilon_{\text{stall}}$, and technical details such as updating the set of active half-spaces \textit{a priori} and throughout the algorithm are omitted.

Since the proposed algorithm achieves the same iterates after escaping the stalling period (lines \ref{alg:start}--\ref{alg:end}) as would be attained with Dykstra's method, the Boyle-Dykstra theorem~\cite{DYKSTRA} is retained and the iterates are guaranteed to converge to the optimal solution $x^\star$ asymptotically. This is formalised in Corollary~\ref{cor:convergence} below.

\begin{cor}[Convergence of Algorithm~\ref{alg:fastforward}~\cite{DYKSTRA}]
\label{cor:convergence}
The sequence of iterates $\{x_m\}$ generated by Algorithm~\ref{alg:fastforward} converges asymptotically to the optimal solution $x^\star = \mathcal{P}_{\mathcal{H}}(x^\circ)$, i.e.\ $\lim_{m\rightarrow\infty}\|x_m-x^\star\|_2=0$.
\end{cor}

\section{Numerical experiment} \label{sec:example}

We demonstrate the efficacy of our method through a numerical experiment. The results of our experiment can be reproduced from the repository at~\cite{DykstraProjectRepo}. The example we used is the projection of an initial point $x^\circ$ on the intersection of a box and a line in $\mathbb{R}^2$. The box is centred at the origin and can be represented as $[-1, 1] \times [-1, 1]$, and the line passes through the points $(0,1)$ and $(2, 0)$. This configuration was first described in detail in~\cite{DYKSTRASTALLING}. The initial point was set to $x^\circ = (-4, 1.4)$ to induce the stalling period, and the ground-truth optimal solution $x^\star$ was computed using the interior-point CQP solver \texttt{quadprog} to serve as a benchmark for error calculations. Half-space activity was monitored by verifying the condition in~\eqref{eq:sets_i} for each set at each iteration, and the convergence behaviour was monitored via the squared $\ell_2$ norm of the absolute distance to optimality, $E(x_m) \eqdef \| x_m - x^\star \|_2^2$.

\definecolor{legreddraw}{HTML}{E41A1C}
\definecolor{legbluedraw}{HTML}{377EB8}
\definecolor{leggreendraw}{HTML}{4DAF4A}
\definecolor{legredfill}{HTML}{B61516}
\definecolor{legbluefill}{HTML}{2C6593}
\definecolor{leggreenfill}{HTML}{3E8C3B}

\begin{figure}[b]
    \centering
\begin{tikzpicture}
  \begin{axis}[
    at={(0,0)},
    anchor=north west,
    width=0.8\linewidth,
    height=4cm,
    xmin=-4.2, xmax=1.2,
    ymin=0.6, ymax=1.5,
    xlabel={$x$ coordinate},
    ylabel={$y$ coordinate},
    grid=both,legend style={fill=white,fill opacity=0.4,text opacity=1,draw=black,cells={anchor=west}},
    legend pos=south west, legend columns=3, transpose legend,
  ]
    \draw[fill=legbluefill, fill opacity=0.2, draw=legbluedraw]
        (axis cs: -1,0.5) rectangle (axis cs: 1,1);
    \addlegendimage{area legend,fill=legbluefill,fill opacity=0.2,draw=legbluedraw}
    \addlegendentry{Box}

    \draw[very thick, legbluedraw] (axis cs:-2,2) -- (axis cs:1,0.5);
    \addlegendimage{line legend,very thick, draw=legbluedraw}
    \addlegendentry{Line}
    
    \draw[line width=2, legredfill!30, line join=round]
      (axis cs:-1,1.4) -- (axis cs:-0.8,1.4) -- (axis cs:-1,1) -- (axis cs:-1,1.4);
    \addlegendimage{line legend,line width=2, draw=legredfill!30}
    \addlegendentry{Stalling}
        
    \addplot[
        black,
        thin, densely dotted, line cap= round,
        mark=*, mark size = 0.75, line join=round,
    ] table[x=x, y=y, col sep=comma] {data_iterates.csv};
    \addlegendentry{Iterates $x_m$}

    \addplot[
    orange,
    only marks,
    mark=star,
    mark size=2.5,
    ] coordinates {(0, 1)};
    \addlegendentry{Projection $x^\star$}

  \end{axis}

  \begin{axis}[
    at={(0,-3.75cm)},
    anchor=north west,
    width=0.8\linewidth,
    height=4cm,
    xmin=0, xmax=30,
    ymode=log,
    ylabel={Squared error},
    xmajorgrids=true,
    xminorgrids=true,
    ymajorgrids=false,
    yminorgrids=false,
    legend pos=south west,
    legend style={fill=white,fill opacity=0.4,text opacity=1,draw=black,cells={anchor=west}},
    legend entries={Original, Algorithm~\ref{alg:fastforward}, Stalling},
  ]
    \draw[line width=2, legredfill!30, line join=round]
      (axis cs:2,0.8) -- (axis cs:16,0.8);

    \addplot[
        black,
        thick,
        mark=*, mark size= 1,
    ] table[x expr=\thisrow{iteration}, y=error, col sep=comma] {data_dykstra_original.csv};
    \addplot[
        leggreendraw,
        thick,
        mark=x, mark size= 2,
    ] table[x expr=\thisrow{iteration}, y=error, col sep=comma] {data_dykstra_modified.csv};
    \addlegendimage{line legend,line width=2, draw=legredfill!30}
  \end{axis}

  \begin{axis}[
    at={(0,-7.5cm)},
    anchor=north west,
    width=0.8\linewidth,
    height=4cm,
    xmin=0, xmax=30,
    ylabel={Half-space activity},
    xlabel={Cycle},
    ytick={0,1},
    grid=both,
    legend style={fill=white,fill opacity=0.4,text opacity=1,draw=black,cells={anchor=west}},
    legend entries={Original, Algorithm~\ref{alg:fastforward}, Stalling},
  ]
    \draw[line width=2, legredfill!30, line join=round]
      (axis cs:2,1) -- (axis cs:16,1);

    \addplot[
        black,
        thick,
        mark=*, mark size= 1,
    ] table[x expr=\thisrow{iteration}, y=halfspace, col sep=comma] {data_dykstra_original.csv};
    \addplot[
        leggreendraw,
        thick,
        mark=x, mark size= 2,
    ] table[x expr=\thisrow{iteration}, y=halfspace, col sep=comma] {data_dykstra_modified.csv};
    \addlegendimage{line legend,line width=2, draw=legredfill!30}
  \end{axis}
\end{tikzpicture}
    \caption{Dykstra's method and Algorithm~\ref{alg:fastforward} applied to the line-box example from~\cite{BAUSCHKESWISS}. Top: The iterates $x_m$ of the two algorithms, with the stalling cycle of Dykstra's method highlighted in red. Middle: Squared errors $E(x_m)$, with the same stalling period highlighted. Bottom: Activity of the vertical half-space (left side of the box), where $1$ indicates that the half-space is active and $0$ indicates that it is inactive; the stalling period is highlighted in red and ends when this half-space becomes inactive.}
    \label{fig:stalling2}
\end{figure}
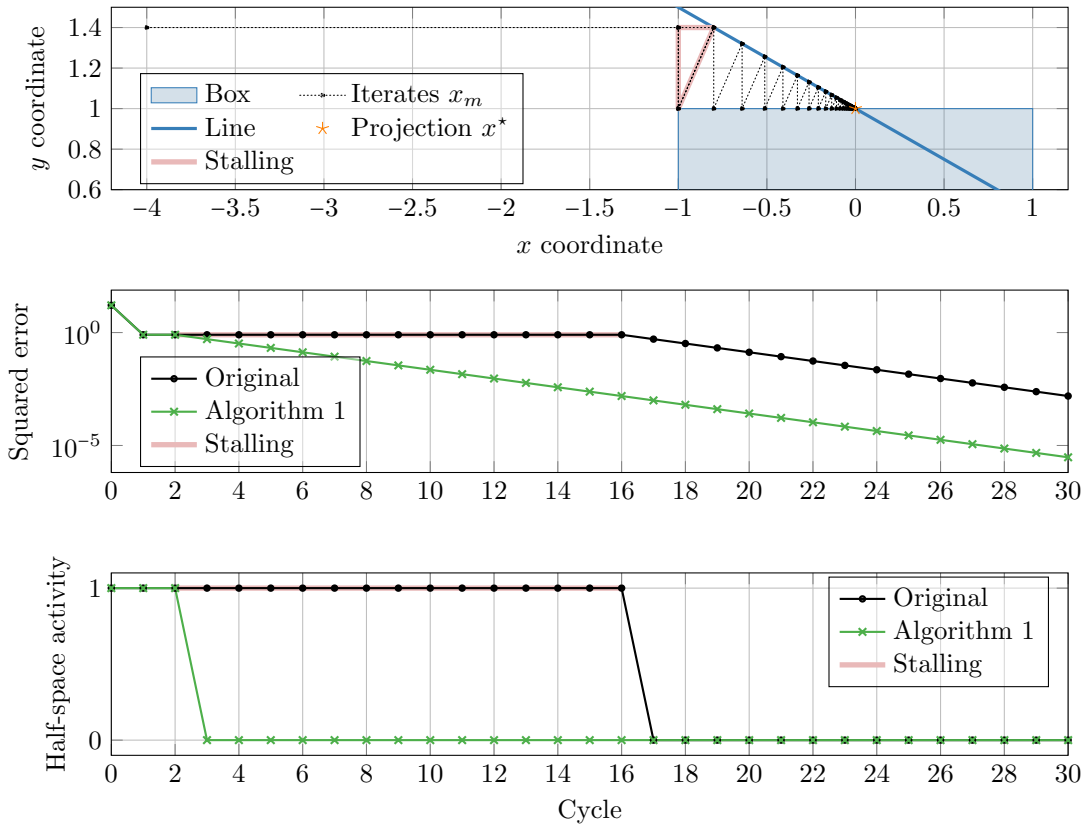

The results for both Dykstra's method and Algorithm~\ref{alg:fastforward} are shown in Figure~\ref{fig:stalling2}. The iterates $x_m$ are shown in the top plot, the errors $E(x_m)$ in the middle plot and the activity of the half-space corresponding to the left side of the box in the bottom plot. Note that the sequences of iterates $x_m$ generated by Dykstra's method and Algorithm~\ref{alg:fastforward} are identical, so they cannot be distinguished in the top plot. However, the error plot (middle) indicates that Dykstra's algorithm stalls up to cycle 16, with the stalled trajectory highlighted in red. The half-space activity shows that stalling ends when the half-space corresponding to the left side of the box becomes inactive.

The middle and bottom plots of Figure~\ref{fig:stalling2} illustrate the application of Algorithm~\ref{alg:fastforward} to the same example, effectively resolving the stalling issue. In Figure~\ref{fig:stalling2}, stalling is detected at cycle 2, upon which the algorithm is fast-forwarded using $N_{stall}=14$ cycles, where $N_{stall}$ is computed using Theorem~\ref{thm:nstall}. At cycle 3, the iterates of the modified algorithm arrive at the same state as the iterates of the original algorithm at cycle 17. After this, the iterates of the modified algorithm become identical to those of the original algorithm. As a result, the modified method exhibits superior convergence properties compared to the standard Dykstra's algorithm. Although this section has showcased a resolution for a simple example in $p = 2$ dimensions with only $n=3$ initially active half-spaces, the results of Theorem~\ref{thm:nstall} and Algorithm~\ref{alg:fastforward} will hold in any number of dimensions and for an arbitrary number of half-spaces.

\section{Conclusion}
\label{sec:conclusion}

Dykstra's algorithm is an iterative method with low per-iteration computational cost for finding Euclidean projections onto the intersection of convex sets, but the stalling phenomenon significantly hinders its practical utility~\cite{BAUSCHKESWISS}. For time-critical applications, such as optimisation in real-time control contexts, an algorithm with an unpredictable and arbitrarily long execution time is not practical. In this paper, we have proposed a solution to the stalling problem for the polyhedral case (Section~\ref{sec:main_result}). We first formalised the stalling period (Definition~\ref{def:stalling}) and then derived its exact length in closed form (Theorem~\ref{thm:nstall}), rendering this quantity computable \textit{a priori} once stalling is detected. Combining this result with the standard Dykstra update, we proposed Algorithm~\ref{alg:fastforward}. This modified algorithm detects stalling and uses the result from Theorem~\ref{thm:nstall} to ``fast-forward'' the algorithm's auxiliary variables past the entire stalling period in a single, computationally inexpensive step. As demonstrated by our numerical results (Section~\ref{sec:example}), the modification eliminates the stalling behaviour (Figure~\ref{fig:stalling2}). Since the proposed algorithm produces the same iterates after escaping the stalling period as would be produced with standard Dykstra, the Boyle-Dykstra theorem is retained and the iterates are guaranteed to converge to the projection $x^\star$ asymptotically (Corollary~\ref{cor:convergence}).
Our method ensures the iterates converge faster, enabling practitioners to implement Dykstra's algorithm in time-critical settings.

We identify several directions for future research. First, extending the analysis from polyhedral sets to broader classes of convex sets would enhance the generality of the results. Second, applying our method to a reformulated linear program in which the optimality conditions are decomposed into a convex and a conic set could provide a benchmark for evaluating the proposed approach. Finally, although we have improved Dykstra's overall convergence properties by eliminating stalling, we have not been able to show that our modified version yields iterates $x_m$ that \emph{monotonically} approach the projection $x^\star$. According to the analysis in~\cite{DYKSTRAPOLY}, the algorithm successively discards sets $\mathcal{H}_i$ for which $x^\star\in\text{int}\,\mathcal{H}_i$, and during these phases, the algorithm can produce iterates that temporarily diverge from $x^\star$. Developing a similar fast-forwarding technique for these non-stalling phases would enable monotonic convergence of the iterates, allow \textit{a priori} bounds on the error $\anynorm{x_m-x^\star}$, and provide deterministic guarantees on solution quality for any fixed iteration budget.

\begin{ack}
This work was supported by Keble College using the Keble College Small Research Grant \#KSRG118.
\end{ack}

\newpage
\printbibliography

\end{document}